\documentclass{amsart}

\usepackage{mathtools}
\usepackage{graphicx} 
\usepackage{amssymb}
\usepackage[hyphens]{url} 
\usepackage{amsmath}
\usepackage{cite}
\usepackage{mathrsfs}
\usepackage{hyperref}
\usepackage{amsthm}
\usepackage{tikz-cd}
\usepackage{caption}
\usepackage{float}
\usepackage[T1]{fontenc}
\usepackage{graphicx}
\usepackage{stmaryrd}
\usepackage{tikz}
\usepackage{tikz-cd} 
\usepackage[all]{xy}
\usepackage{comment}
\usepackage{bm}
\usepackage{comment}
\usepackage{amsmath,amsfonts,amssymb}
\usepackage{tcolorbox}
\usepackage{cleveref}
\usepackage{hyperref}
\usepackage[a4paper]{geometry}
\usepackage{nameref}
\usepackage{mathabx}
\usepackage{enumitem}
\usepackage{xcolor}
\usepackage{nicefrac}
\usepackage{array}
\hypersetup{
	colorlinks=true,
	linkcolor={blue},
	citecolor={blue},
	urlcolor={blue}}

\theoremstyle{plain}
\newtheorem*{mainthm}{Main Theorem}

\newtheorem*{conj*}{Conjecture}

\newtheorem*{cor*}{Corollary}

\theoremstyle{definition}

\newtheorem*{ex*}{Example}

\theoremstyle{remark}

\newtheorem*{rmk*}{Remark}

\newlength{\plarg}
\usetikzlibrary{cd}

\newcommand{\abs}[1]{\left|#1\right|}
\newcommand{\norm}[1]{\left\lVert#1\right\rVert}
\newcommand{\ab}{\mathrm{ab}}

\numberwithin{equation}{section}

\title[Virtual homological torsion: abundance versus growth]{Virtual homological torsion: abundance versus growth in books of $I$-bundles}
\author{Jonathan Fruchter}

\begin{document}

\begin{abstract}
    Let $\mathcal{B}$ be a book of $I$-bundles, all of whose pages are surfaces of negative Euler characteristic. In this short note, we prove that torsion in the first homology of $\mathcal{B}$ grows subexponentially in the index along any exhausting tower of regular finite-sheeted covers. By contrast, recent work of Ascari and the author shows that, apart from the obvious exceptions, $\mathcal{B}$ has abundant virtual homological torsion, which can grow exponentially along exhausting towers of non-regular finite covers.
\end{abstract}

\maketitle

The behaviour of torsion in the homology of finite-sheeted covers of manifolds has emerged as a central topic in geometric topology. In particular, within the world of $3$-manifolds, homological torsion growth has become an important theme in the interaction of geometry, topology, and number theory. Influential work of Bergeron--Venkatesh and L\"uck has drawn attention to the asymptotic growth of torsion along exhausting towers of normal finite-sheeted covers. They conjecture:

\begin{conj*}[\hspace{-.05em}{\cite{BV}, \cite{Lck1}, \cite{Lck2}}]
    Let $M$ be a compact connected orientable and irreducible $3$-manifold with an infinite fundamental group and empty or toroidal boundary. Then $M$ admits a cofinal tower of regular finite-sheeted covers
    \[M=M_0 \twoheadleftarrow M_1 \twoheadleftarrow M_2 \twoheadleftarrow \cdots \]
    for which
    \[
    \lim_{n \rightarrow \infty } \frac{\log(\abs{\mathrm{Tor}(H_1(M_n;\mathbb{Z}))})}{[\pi_1(M):\pi_1(M_n)]}=\frac{\sum_{i=1}^n\mathrm{vol}(P_i)}{6\pi},
    \]
    where the sum ranges over the hyperbolic pieces $P_i$ in the JSJ decomposition of $M$.
\end{conj*}

We remark that this is a special case of L\"uck's approximation conjecture for $L^2$-torsion \cite[Conjecture 1.11]{Lck2}. While L\"uck's conjecture has been confirmed in certain settings, no finitely presented residually finite group is currently known to exhibit exponential homological torsion growth in dimension one along a residual chain of finite-index normal subgroups. \par \smallskip

On the other hand, several results show a rich supply of virtual homological torsion. A space (or a group) is said to have \emph{abundant virtual homological torsion} if every finite abelian group appears as a direct factor in the first homology of one of its finite covering spaces (or finite-index subgroups). Sun \cite[Theorem 1.5]{Sun} and Chu--Groves \cite[Theorem 1.1]{Chu} proved that any non-graph $3$-manifold with empty or toroidal boundary possesses this property. More recently Ascari and the author \cite[Theorem A]{us} established the same for hyperbolic groups that split as graphs of free groups with cyclic edge groups, other than free products of free and surface groups. \par \smallskip

The aim of this note is to highlight the tension between abundance and growth of virtual homological torsion by considering \emph{books of $I$-bundles}, introduced by Culler and Shalen. Following \cite{books}, a book of $I$-bundles $\mathcal{B}$ is a compact, connected and orientable $3$-manifold constructed from a finite collection of solid tori (the \emph{bindings}) and a finite collection of $I$-bundles over connected compact surfaces with boundary (the \emph{pages}). For each page $P$, which is an $I$-bundle over a surface $\Sigma$, every \emph{vertical boundary annulus} in the preimage of $\partial \Sigma$ in $P$ is attached to a chosen homotopically essential annulus in the boundary of some binding $B$. The group $\pi_1(\mathcal{B})$ is finitely presented, residually finite, and contains no infinite normal amenable subgroups \cite[Corollary 0.6]{Gromov}. With few exceptions, the aforementioned results of Ascari and the author imply that books of $I$-bundles also have abundant virtual homological torsion. Nevertheless, we prove:

\begin{mainthm} \label{mainthm}
    Let $\mathcal{B}$ be a book of $I$-bundles, all of whose pages are surfaces of negative Euler characteristic. If $\mathcal{B}=\mathcal{B}_0\twoheadleftarrow\mathcal{B}_1\twoheadleftarrow \cdots $ is any tower of regular finite-sheeted covers such that $\bigcap_n \pi_1(\mathcal{B}_n)=1$, then
    \[
    \lim_{n \rightarrow \infty } \frac{\log(\abs{\mathrm{Tor}(H_1(\mathcal{B}_n;\mathbb{Z}))})}{[\pi_1(\mathcal{B}):\pi_1(\mathcal{B}_n)]}=0.
    \]
\end{mainthm}

\begin{rmk*}
    It follows from the work of Bridson and Kochloukova \cite[Corollary C]{BK} that the first $L^2$-Betti number of such $\mathcal{B}$ never vanishes, and therefore its $L^2$-torsion is undefined.
\end{rmk*}

Two particularly fruitful strategies for proving subexponential homological torsion growth are L\"{u}ck’s theorem \cite[Theorem 1.14]{Lck2}, which applies to a large class of finitely presented residually finite groups with infinite amenable normal subgroups, and the \emph{effective rebuilding} method of Abert--Bergeron--Fr\k{a}czyk--Gaboriau \cite{rebuild}. Books of $I$-bundles are impervious to both approaches, and our proof therefore uses a simpler more direct argument. \par \smallskip

Indeed, $\pi_1(\mathcal{B})$ decomposes as a graph of groups with free vertex groups and cyclic edge groups, and it is hyperbolic \cite{BF}. This places it outside the reach of the methods from the previous paragraph (see, in particular, \cite[Theorem~H and the subsequent remark]{rebuild}), yet within the scope of \cite[Theorem A]{us}. Hence, unless $\pi_1(\mathcal{B})$ is a free product of free and surface groups, it has abundant virtual homological torsion; in fact, torsion can even grow exponentially in the index along carefully chosen exhausting sequences of non-normal finite-index subgroups \cite[Main Theorem]{us2}. \par \smallskip

Taken together, these results highlight an interesting tension: while an abundance of torsion, and even exponential growth, are guaranteed in sufficiently flexible families of covers, the behaviour in normal towers is much more delicate, and homological torsion fails to grow at an exponential rate. In the conjectural landscape of Bergeron--Venkatesh, such behaviour is reserved for higher–dimensional arithmetic groups: for instance, finite–volume arithmetic hyperbolic manifolds of simplest type in dimension $n>3$ have abundant virtual homological torsion \cite[Corollary~1.2]{Chu}, which is conjectured to grow only subexponentially \cite{BV}. In sharp contrast, books of $I$–bundles are essentially $2$–dimensional objects in spirit. \par \smallskip

We are not aware of other finitely presented, residually finite, one–ended groups without infinite amenable normal subgroups where these phenomena coexist. If one allows for infinite normal amenable subgroups, which under mild finiteness assumptions forces subexponential torsion growth \cite[Theorem 1.14]{Lck2}, degenerate examples appear (e.g. $\pi_1(M)\times \mathbb{Z}$, where $M$ is a non–graph $3$–manifold with empty or toroidal boundary). More robust examples with \emph{unbounded} virtual homological torsion are also not difficult to construct. For example, let $A\in \mathrm{GL}(2,\mathbb{Z})$ be a matrix with two real eigenvalues $\lambda>1$ and $\lambda^{-1}$, and let $G$ be the metabelian group $\mathbb{Z}^2 \rtimes_A \langle t \rangle$. Unwrapping the monodromy $A$ gives a (non-exhausting) sequence of cyclic covers $G_n$ of $G$ in which $\abs{\mathrm{Tor}(G_n^{\ab})}$ diverges as $n\to\infty$. Note, however, that the virtual homological torsion of $G$ is always $2$-generated, and thus $G$ does not have abundant virtual homological torsion.

\begin{rmk*}
    Another potential source of examples are polynomially growing free-by-cyclic groups. These have the cheap $\alpha$-rebuilding property for every $\alpha$, and therefore admit subexponential homological torsion growth \cite{sam, sam2}, but it is not known whether they have abundant virtual homological torsion\footnote{This question was raised by Cheetham-West, and appears in the American Institute of Mathematics' problem list ``\href{http://aimpl.org/freebycyclic/}{Rigidity properties of free-by-cyclic groups}''.}
\end{rmk*}

\subsection*{Proof setup} In preparation for the proof, denote the bindings of $\mathcal{B}$ by $B_1,\ldots, \allowbreak B_m$, and its pages by $P_1,\ldots, P_k$ with underlying surfaces $\Sigma_1,\ldots, \Sigma_k$. As in \cite[p.~286]{books}, contracting each $B_i$ onto its longitude, and each $P_i$ onto $\Sigma_i$, we obtain a $2$-dimensional CW-complex $X$ with $\pi_1(X)=\pi_1(\mathcal{B})$. \par \smallskip

The space $X$ admits a \emph{graph of spaces decomposition} $\mathcal{X}$ (in the sense of Scott and Wall \cite{Scott}), whose vertex spaces are either circles $C_1,\ldots,C_m$ coming from the bindings of $\mathcal{B}$, or compact surfaces with boundary of negative Euler characteristic $\Sigma_1,\ldots,\Sigma_k$ coming from its pages. The edge spaces of $\mathcal{X}$ are all circles, and the gluings in $\mathcal{B}$ dictate the underlying graph $\Gamma$ of $\mathcal{X}$ as well as its edge maps, which either identify an edge space $X_e$ with a boundary component of some $\Sigma_i$, or attach $X_e$ to some $C_j$ by a covering map of degree at least $1$. This decomposition of $X$ passes naturally to its covering spaces, endowing them with a similar combinatorial description that we will exploit in our proof.\par \smallskip

Indeed, consider a covering space $p:\widehat{X}\twoheadrightarrow X$. The graph of spaces decomposition $\mathcal{X}$ lifts to a decomposition $\widehat{\mathcal{X}}$ of $\widehat{X}$. Its vertex and edge spaces are the connected components of the $p$-preimages of the vertex and edge spaces of $\mathcal{X}$. Its edge maps are given by the \emph{elevations} of the edge maps of $\mathcal{X}$ to $\widehat{\mathcal{X}}$ (see Wise's work for further detail {\cite[Definition 2.18]{wise}}). Briefly, the elevations of an edge map $i_e:X_e \longrightarrow X_v$ to $\widehat{X}$ are given by the restrictions $\widehat{i}_j$ of the map $\widehat{i}_e$ to the connected components $\widehat{X}_{e_j}$ of $X_e \times_X \widehat{X}$ arising from the pullback diagram
\begin{center}
\begin{tikzpicture}
    \node (A) at (0,1.3) {$X_e \times_X \widehat{X}=\bigsqcup_{j=1}^m \widehat{X}_{e_j}$};
    \node (B) at (5,1.3) {$\widehat{X}$};
    \node (C) at (0,0) {$X_e$};
    \node (D) at (5,0) {$X$};

    \draw[->] (A) -- node[above] {$\widehat{i}_e=\bigsqcup_{j=1}^m \widehat{i}_j$} (B);
    \draw[->] (B) -- node[right] {$p$} (D);
    \draw[->] (A) -- node[left] {$q=\bigsqcup_{j=1}^m$} (C);
    \draw[->] (C) -- node[below] {$i_e$} (D);
\end{tikzpicture}
\end{center}
Finally, suppose the covering map $\widehat{X}\twoheadrightarrow X$ has finite degree. Let $\widehat{C}_i$ be a lift of a circle vertex space $C_i$ to $\widehat{X}$, and let $\widehat{i}_1,\ldots,\widehat{i}_\ell$ be the elevations of an edge map $i_e:X_e \longrightarrow C_i$ with targets in $\widehat{C}_i$. Each $\widehat{i}_j$ is a covering map whose degree divides that of $i_e$ and
\[
\deg(\widehat{i}_1)+\cdots + \deg(\widehat{i}_\ell)=\deg(i_e).
\]
Similarly, if $\widehat{\Sigma}_j$ lies above a surface vertex space $\Sigma_j$, then the elevations with target in $\widehat{\Sigma}_j$ identify each component of $\partial \widehat{\Sigma}_j$ with a unique edge space of $\widehat{\mathcal{X}}$. \par \smallskip
With this setup in hand, we turn to the proof of the main theorem. \par \smallskip

\begin{proof}[Proof of the {\hyperref[mainthm]{Main Theorem}}]
    Let 
    \[
    \mathcal{B}=\mathcal{B}_0\twoheadleftarrow \mathcal{B}_1 \twoheadleftarrow\mathcal{B}_2\twoheadleftarrow \cdots
    \]
    be a cofinal tower of regular finite-sheeted covers of $\mathcal{B}$, and let
    \[
    X=X_0 \twoheadleftarrow X_1 \twoheadleftarrow X_2 \twoheadleftarrow \cdots
    \]
    be the corresponding tower of regular finite-sheeted covers of $X$. \par \smallskip

    For each $n$, the degree of the covering map $X_n \twoheadrightarrow X$ splits into two contributions: one coming from the \emph{configuration} of $\mathcal{X}_n$ (that is, the complexity of the underlying graph $\Gamma_n$), and another coming from \emph{unwrapping} the various vertex and edge spaces of $\mathcal{X}$. Since the tower $X\twoheadleftarrow X_1 \twoheadleftarrow X_2 \twoheadleftarrow\cdots$ is cofinal, both contributions diverge as $n\longrightarrow \infty$. We will show that while torsion in $H_1(X_n;\mathbb{Z})$ can potentially grow exponentially with respect to the configuration part of the index, once the unwrapping contribution is taken into account the overall growth is forced to be subexponential. \par \smallskip
    
    We begin by fixing necessary notation. For every circle vertex space $C_i$ of $\mathcal{X}$, fix a generator $t_{i}$ of $\pi_1(C_i)\cong \mathbb{Z}$, and let $\mathrm{val}_i$ denote its valence in the underlying graph $\Gamma$ of $\mathcal{X}$. For every edge $e\in \mathrm{E}(\Gamma)$ that terminates at a circle vertex space $C_i$, write $d_e$ for the degree of the covering map $i_e:X_e \twoheadrightarrow C_i$. Further, for each surface vertex space $\Sigma_j$ of $\mathcal{X}$, write $\sigma_{j,1},\ldots, \sigma_{j,s_j}$ for the connected components of $\partial \Sigma_j$, each of which is endowed with an orientation. \par \smallskip
    
    Consider a covering space $X_n$ in the sequence. For every circle vertex space $C_i$ of $\mathcal{X}$, let $C^n_{i,1},\ldots,C^n_{i,\ell _i}$ be the circle vertex spaces of $\mathcal{X}_n$ that cover $C_i$. Likewise, let $\Sigma^n_{j,1},\ldots, \Sigma^n_{j,p_j}$ denote the lifts of the surface vertex space $\Sigma_j$ of $\mathcal{X}$ to $\mathcal{X}_n$. We also fix a generator $t^n_{i,j}$ for every circle vertex space $C^n_{i,j}$ of $\mathcal{X}_n$, and denote the boundary of each surface vertex space by
    \[
    \partial \Sigma^n_{j,i} = \left( \sigma^n_{j,i,k}\right)_{k=1}^{s_{j}},
    \]
    with orientations inherited from $\Sigma_j$. Let 
    \[(x^n_{j,i,k})_{k=1}^{r_j}\]
    be the standard handle or crosscap surface generators, so that together with $( \sigma^n_{j,i,k})_{k=1}^{s_{j}}$ they generate $\pi_1(\Sigma^n_{j,i})$. Since the cover $X_n\twoheadrightarrow X$ is regular, the surfaces $\Sigma^n_{j,1},\ldots, \Sigma^n_{j,p_j}$ are all homeomorphic, and $s_j, r_j$ above do not depend on the chosen lift $\Sigma^n_{j,i}$ of $\Sigma_j$. \par \smallskip

    Note that $\ell _i$, $s_j$ and $r_j$ depend on $n$. To keep notation light we omit this dependence from our notation, but when needed, we will write $\ell_i=\ell_i(n)$. \par \smallskip
    
    By Mayer-Vietoris, the first homology of $X_n$ is given by subjecting the free abelian group
    \[
    \bigoplus_{j=1}^{k} \bigoplus_{i=1}^{p_j} H_1(\Sigma^n_{j,i}) \; \oplus \; \langle t^n_{i,j} \; \vert \; i\le m,\;j\le \ell_i\rangle
    \]
    to the following relations:
    \begin{enumerate}
        \item for every edge $e \in \Gamma_n$ that attaches the boundary component $\sigma^n_{j,i,k}$ of $\Sigma^n_{j,i}$ to a circle $C^n_{i',j'}$, a relation
        \begin{equation}\label{rel1}
        \sigma^n_{j,i,k}=d_e \cdot t^n_{i',j'}
        \end{equation}
        where $d_e$ is the signed degree of the attaching map $i_e:(X_n)_e\twoheadrightarrow C^n_{i',j'}$. Note that if $e$ lies above the edge $e'$ of $\Gamma$, then $d_e \mid d_{e'}$. 
        \item for every orientable surface vertex space $\Sigma^n_{j,i}$ of $\mathcal{X}_n$, a relation
        \begin{equation} \label{rel2}
        \pm \sigma^n_{j,i,1} \pm \cdots \pm \sigma^n_{j,i,s_j}=0,
        \end{equation}
        where the signs are determined by the previously chosen orientations. We remark that if $\Sigma^n_{j,i,1}$ is a punctured sphere, we view its fundamental group as generated by its boundary components with the relation that their product is trivial.
        \item for every non-orientable surface vertex space $\Sigma^n_{j,i}$ of $\mathcal{X}_n$, a relation
        \begin{equation} \label{rel3}
            2\cdot x^n_{j,i,1} + \cdots + 2\cdot x^n_{j,i,r_j} \pm \sigma^n_{j,i,1} \pm \cdots \pm \sigma^n_{j,i,s_j}=0,
        \end{equation}
        where again, the signs in the sum are determined by the chosen orientations.
    \end{enumerate}

    By the definition of a book of $I$-bundles, every boundary component of each surface vertex of $\mathcal{X}$, and hence of $\mathcal{X}_n$, participates in the gluing, and is attached to a unique circle vertex space. This implies that 
    \[H_1(X_n)=A_n \oplus \mathbb{Z}^{r_n},\]
    where $A_n$ is a finitely generated abelian group with presentation obtained by plugging the relations in \ref{rel1} into those in \ref{rel2} and \ref{rel3}. In more detail, the group $A_n$ is described by the matrix $\mathcal{A}_n$ constructed as follows:
    \begin{enumerate}
        \item The rows of $\mathcal{A}_n$ correspond to the surface vertex spaces $\Sigma^n_{j,i}$ of $\mathcal{X}_n$,
        \item The columns of $\mathcal{A}_n$ are given by
        \begin{enumerate}
            \item the circle vertex spaces $C^n_{i',j'}$ of $\mathcal{X}_n$, and
            \item the crosscap generators $(x^n_{j,i,k})_{k=1}^{r_j}$ of every non-orientable surface vertex space $\Sigma^n_{j,i}$ of $\mathcal{X}_n$.
        \end{enumerate} 
        \item Each row encodes the boundary relations \ref{rel2} and \ref{rel3} of a surface vertex $\Sigma^n_{j,i}$ by replacing each $\sigma^n_{j,i,k}$ with $d_e\cdot t^n_{i',j'}$ if it is attached to the circle vertex $C^n_{i',j'}$ by an edge $e$.
    \end{enumerate}

    By Smith normal form, the size of the torsion subgroup of $A_n$ is the product of the invariant factors of $\mathcal{A}_n$. This in turn, is given by the greatest common divisor of the determinants of all full-rank $\mathrm{rk}(\mathcal{A}_n) \times \mathrm{rk}(\mathcal{A}_n)$ minors of $\mathcal{A}_n$. Thus, bounding $\abs{\mathrm{Tor}(H_1(X_n))}=\abs{\mathrm{Tor}(A_n)}$ reduces to bounding the determinants of all such minors. To this end, we appeal to Hadamard's inequality, which bounds the determinant of a matrix by the product of the lengths of its column vectors. \par \smallskip
    For each column $C^n_{i,j}$ of $\mathcal{A}_n$ we have that
    \begin{align*}
    \norm{C^n_{i,j}}_2 & \le \norm{C^n_{i,j}}_1 \\
    &= \sum_{\substack{e\;\text{adjacent} \\ \text{to }C^n_{i,j}\;\text{in }\Gamma_n}} \abs{d_e} \\
    & = \sum_{\substack{e'\;\text{adjacent} \\ \text{to }C_i\;\text{in }\Gamma}} \abs{d_{e'}} \\
    &\le \max_{i=1}^m(\mathrm{val}_i)\cdot \max\{\abs{d_{e'}}\;\vert\;e'\in \mathrm{E}(\Gamma)\},
    \end{align*}
    and since each column $x^n_{j,i,k}$ given by a crosscap generator of a non-orientable surface has a single nonzero entry (which is equal to $2$) at the row $\Sigma^n_{j,i}$,
    \[
    \norm{x^n_{j,i,k}}_2 =2.
    \]
    The number of columns of $\mathcal{A}_n$ is
    \[
    \sum_{i=1}^m \ell_i + \sum_{\Sigma^n_{i,j}\;\text{non orientable}} r_j,
    \]
    whereas its number of rows is bounded from above by
    \begin{align*}
    \mathrm{rows}(\mathcal{A}_n) &\le \sum_{i=1}^m \ell_i \cdot \max\{\abs{d_{e}}\;\vert\;e \text{ is adjacent to }C_i \text{ in } \Gamma\} \\
    & \le m\cdot \max\{\ell_1,\ldots,\ell_n\}\cdot \max\{\abs{d_{e'}}\;\vert\;e'\in \mathrm{E}(\Gamma)\}.
    \end{align*}
    The inequality above follows from the fact that if $C_i$ and $\Sigma_j$ were connected by an edge $e$ in $\mathcal{X}$, then every lift $C^n_{i,j}$ of $C_i$ is adjacent to at most $d_e$ lifts of $\Sigma_j$ to $\mathcal{X}_n$. Write
    \[\mathrm{val}=\max_{i=1}^m(\mathrm{val}_i)\;\;\text{and}\;\;d=\max\{\abs{d_{e'}}\;\vert\;e'\in \mathrm{E}(\Gamma)\}.\]
    We may now apply Hadamard's inequality using the bound on the number of rows as a bound for the rank of $\mathcal{A}_n$:
    \begin{align*}
        \abs{\mathrm{Tor}(H_1(X_n;\mathbb{Z}))} &=\abs{\mathrm{Tor}(A_n)}\\
        &=\mathrm{gcd}\{\det(\mathcal{A})\;\vert\;\mathcal{A}\;\text{is a full-rank } \mathrm{rk}(\mathcal{A}_n)\times \mathrm{rk}(\mathcal{A}_n) \text{ minor of } \mathcal{A}_n\} \\
        & \le d\cdot \prod_{i=1}^m \prod_{j=1}^{\ell_i} \max\{\norm{C^n_{i,j}}_2, 2\} \\
        & \le (2 \mathrm{val} \cdot d)^{d\cdot m\cdot \max\{\ell_1,\ldots,\ell_n\}}.
    \end{align*}
    
    Recall that each $\ell_i=\ell_i(n)$ depends on $n$. For every $n$, let $i_n\le m$ be such that 
    \[\ell_{i_n} = \ell_{i_n}(n)= \max\{\ell_1(n),\ldots,\ell_m(n)\}.\]
    Let $D_n$ be the degree of the covering map $C^n_{{i_n},j}\twoheadrightarrow C_{i_n}$ (which does not depend on $j$ since $X_n$ is a regular cover of $X$). We have that
    \[
    [\pi_1(X):\pi_1(X_n)]=\ell_{i_n} \cdot D_n,
    \]
    and since our sequence of covers is cofinal, $\lim_{n\rightarrow \infty}D_n = \infty$. We remark that here, $\ell_{i_n}$ and $D_n$ represent, respectively, the configuration and unwrapping components of the index mentioned earlier. Putting everything together,
    \begin{align*}
        \lim_{n \rightarrow \infty } \frac{\log(\abs{\mathrm{Tor}(H_1(\mathcal{B}_n;\mathbb{Z}))})}{[\pi_1(\mathcal{B}):\pi_1(\mathcal{B}_n)]} & = \lim_{n\rightarrow \infty}\frac{\log\mathrm{Tor}(A_n)}{[\pi_1(X):\pi_1(X_n)]} \\
        & \le \lim_{n \rightarrow \infty } \frac{\log\left((2 \mathrm{val}\cdot d)^{d\cdot m\cdot \ell_{i_n}}\right)}{\ell_{i_n}\cdot D_n} \\
        & =  \lim_{n \rightarrow \infty }\frac{\log(2\mathrm{val}\cdot d)\cdot d\cdot m \cdot \ell_{i_n}}{\ell_{i_n}\cdot D_n}\\
        & =  \lim_{n \rightarrow \infty }\frac{\log(2\mathrm{val}\cdot d)\cdot dm}{D_n}=0,
    \end{align*} \end{proof}

\subsection*{Acknowledgements} This work received funding from the European Union (ERC, SATURN, 101076148) and the Deutsche Forschungsgemeinschaft (DFG, German Research Foundation) under Germany's Excellence Strategy - EXC-2047/1 - 390685813. I warmly thank Cameron Gates-Rudd, Sam Hughes, Lawk Mineh, Andrew Ng and Nikolay Nikolov for their comments and suggestions.

\bibliographystyle{plain}

\vspace{.5cm}

\textsc{Mathematisches Institut, Rheinische Friedrich-Wilhelms-Universit\"at Bonn, Endenicher Allee 60, 53115 Bonn, Germany}

\emph{Email address:} \texttt{fruchter@math.uni-bonn.de}
\end{document}